\documentclass[10pt]{article}
 \usepackage{amssymb}

\def\endofproof{\hfill\vrule height6pt width5pt depth0pt}
\def\Proof{{\bf Proof :}~}
\def\bC{{\mathbb C}} \def\bF{{\mathbb F}}  
\def\bQ{{\mathbb Q}}    
  \def\H{{\cal H}} 
\def\bfC{{\cal C}}

\def\t{\tilde}

\def\C{C(X,a,b)}
\def\p{\pi_1 (X,a)}

\def\E{\mathbb E}

\def\P{{\mathbb P}^1}  
\def\M{{\cal M}}

\newtheorem {Theorem}{Theorem}[section]

\newtheorem {Definition}[Theorem]{Definition}

\newtheorem {Remark}[Theorem]{Remark}
\newtheorem {Corollary}[Theorem]{Corollary}
\newtheorem {Lemma}[Theorem]{Lemma}

\newcommand{\beq}{\begin{equation}}
\newcommand{\eeq}{\end{equation}}

\begin{document}
\setcounter{section}{0}

 \def\fge{{$L/k(x)$\ }}   \def\s{{\sigma}}
 \def\k{\kappa} \def\L{\Lambda}  \def\C{{\bf C}}
  \def\bF{{\mathbb F}}  \def\qed{\endofproof}
  \def\bt{{\bf t}}   \def\bfC{{\bf C}}
  \def\HC{{\cal H}_r(\C)}  \def\E{{\cal E}} \def\S{{\cal S}}
  \def\SA{{\cal S}_r^{(\A)}(G)} \def\orr{{\cal O}(r+1)}
  \def\TG{{\cal T}_r(G)} \def\oor{{\cal O}_r} \def\T{{\cal T}}
   \def\K{{\cal K}} \def\k{\kappa}
    \def\A{{\bf A}}  \def\UR{{\cal U}^{(r)}}
    \def\Ni{{\mbox{Ni}}}  \def\Vtt{(\tilde V)\tilde{}}
     \def\Vht{(\hat V)\tilde{}}
     \def\EC{{\cal E}(\C)}



     \centerline{\Large The monodromy group of a function on a
     general curve}
     \bigskip

     \centerline{Kay Magaard
     \footnote{Partially supported by NSA grant MDA-9049810020} and Helmut V\"olklein
     \footnote{Partially supported by NSF grant DMS-0200225}}
 \bigskip

  \centerline{Wayne State University and University of Florida }

 \bigskip

\centerline{October 2001}

     \vspace{2cm}
     \noindent {\bf Abstract:} 
{\it Let $C_g$ be a general curve of genus $g\ge4$. Guralnick and
others proved that the monodromy group of a cover $C_g\to\P$
of degree $n$ is either $S_n$ or $A_n$. We show that $A_n$
occurs for $n\ge
2g+1$. The corresponding result for $S_n$ is classical. }

\def\p{{\bf p}}
 \def\BB{{\cal B}} \def\FF{{\cal F}}
 \def\Hs{{\cal H}_\s} \def\Hss{{\cal H}_{\ss}} \def\EE{{\cal E}} \def\NN{{\cal N}}
\def\HHR{{\cal H}^{(r)}} \def\HHr{{\cal H}_{r}} \def\H{{\cal H}}
\def\U{{\cal U}}  \def\UR{{\cal U}^{(r)}} \def\Ur{{\cal U}_{r}}
\def\md{\mbox{{\bf mod-dim}}}  \def\ss{\tilde\s}

\section{Introduction}

\bigskip
\noindent 
Let $C_g$ be a general curve of genus $g\ge2$ (over $\bC$). Then
$C_g$ has a cover to $\P$ of degree $n$ if and only if $2(n-1)\ge
g$. This is a classical fact of algebraic geometry. (It is part of
Brill-Noether theory, which more generally considers maps of a
curve to ${\mathbb P}^m$, see \cite{HM}, Ch. 5). If $C_g$ has a cover
to $\P$ of degree $n$, then there is such a cover that is  simple,
i.e., has monodromy group $S_n$ and all inertia groups  are
generated by transpositions. The question arises whether $C_g$
admits other types of covers to $\P$.

If there is a cover $C_g\to\P$ branched at $r$ points of $\P$ and
$g\ge2$ then $r\ge 3g$ (see Remark \ref{nec} below). Zariski
\cite{Z} used this to show that if $g>6$ then there is no such
cover with solvable monodromy group. He made a conjecture on the
existence of such covers for $g\le6$, but there is a
counterexample to that, see Fried \cite{Fr}, Fried/Guralnick
\cite{FrGu}.

The condition $r\ge 3g$ was further used by Guralnick to restrict
the possibilities for the monodromy group $G$ of a cover
$C_g\to\P$ of degree $n$. Assume the cover does not factor
non-trivially, i.e.,  $G$ is a primitive subgroup of $S_n$.
(Knowledge of this case is sufficient to know all types of covers
$C_g\to\P$; this was already observed by Zariski \cite{Z}, see
\cite{GM}). If further $g>3$, then $G=S_n$ or $G=A_n$.
For $g=3$ there are 3 additional cases, with $n=7,8,16$
and $G=GL_3(2), AGL_3(2),AGL_4(2)$, respectively. This was proved
by Guralnick and Magaard \cite{GM} and Guralnick and Shareshian \cite{GS},
using the classification of
finite simple groups. 
There is also
a corresponding result for $g=2$, but it is less definitive.

As noted in \cite{GM}, it was not known whether the case $G=A_n$
actually occurs. This is answered in the affirmative in this paper.
More precisely, we prove the following: Let $g \geq 3$ and $n \ge 2$. Then
the general curve of genus $g$ admits a cover to $\P$
of degree $n$ with monodromy group $A_n$ such that all inertia
groups are generated by double transpositions
 if and only if $n\ge
2g+1$. The same statement holds when we replace double transpositions 
by $3$-cycles (see Theorem \ref{thm}). We refine the latter result in
Theorem \ref{thm3} by showing that both of the two types of $3$-cycle
covers occur for the general curve. (See Fried \cite{Fr1} and
Serre \cite{Se1}, \cite{Se2} for this type distinction).
We also study the exceptional cases in genus 3. 

A preliminary version of this paper has been circulated since
October 2001. It was brought to our attention that in 
a recent preprint S. Schr\"oer \cite{Schr} proves  
a weaker version of our result on $3$-cycles (which, however,
also holds in positive characteristic): The 
locus in $\M_g$ of curves admitting a cover to 
$\P$ with only triple ramification points has dimension
$\geq max(2g-3,g)$.  
 
The authors are grateful to Bob Guralnick for raising the
question of moduli dimension for alternating groups. 
We further gratefully acknowledge the kind support of Gerhard Frey without
whose expert advice this paper may have never been finished.

\section{Moduli dimension of a tuple in $S_n$}

\bigskip
\noindent Let $\P=\P_\bC$ the Riemann sphere. Let $\UR$ be the
open subvariety of $(\P)^r$ consisting of all $(p_1,...,p_r)$ with
$p_i\ne p_j$ for $i\ne j$. Consider a cover $f:X\to\P$ of degree
$n$, with branch points  $p_1,...,p_r\in\P$. Pick $p\in
\P\setminus\{p_1,...,p_r\}$, and choose loops $\gamma_i$ around
$p_i$ such that $\gamma_1,...,\gamma_r$ is a standard generating
system of the fundamental group $\Gamma:=\pi_1(
\P\setminus\{p_1,...,p_r\},p)$ (see \cite{Buch}, Thm. 4.27); in
particular, we have $\gamma_1\cdots\gamma_r=1$. Such a system
$\gamma_1,...,\gamma_r$ is called a homotopy basis of
$\P\setminus\{p_1,...,p_r\}$. The group $\Gamma$ acts on the fiber
$f^{-1}(p)$ by path lifting, inducing a transitive subgroup $G$ of
the symmetric group $S_n$ (determined by $f$ up to conjugacy in
$S_n$). It is called the {\bf monodromy group} of $f$. The images
of $\gamma_1,...,\gamma_r$ in $S_n$ form a tuple of permutations
called a tuple of {\bf branch cycles} of $f$.

Let $\s_1,...,\s_r$ be elements $\ne1$ of  the symmetric group
$S_n$ with $\s_1\cdots\s_r=1$, generating a transitive subgroup.
 Let $\s=(\s_1,...,\s_r)$. We call such a tuple {\bf admissible}.
We say a cover $f:X\to\P$ of degree $n$ is of type $\s$ if it has
$\s$ as tuple of branch cycles relative to some homotopy basis of
$\P$ minus the branch points of $f$. The genus $g$ of $X$ depends
only on $\s$ (by the Riemann-Hurwitz formula); we write $g=g_\s$.

Let $\Hs$ be the set of pairs $([f], (p_1,...,p_r))$, where $[f]$
is an equivalence class of covers of type $\s$, and $p_1,...,p_r$
is an ordering of the branch points of $f$. We use the usual notion of 
equivalence of covers, see \cite{Buch}, p. 67. Let $\Psi_\s:\Hs\to
\UR$ be the map forgetting $[f]$. The {\bf Hurwitz space} $\Hs$
carries a natural structure of quasiprojective variety such that
$\Psi_\s$ is an algebraic morphism, and an unramified covering in the
complex topology (see \cite{FrV},\cite{Buch}, \cite{We}). We also have the
morphism $$\Phi_\s:\ \Hs\ \ \to\ \ \M_g$$ mapping $([f],
(p_1,...,p_r))$ to the class of $X$ in the moduli space $\M_g$
(where $g=g_\s$). Each irreducible component of $\Hs$ has the same
 image in
$\M_g$ (since the action of $S_r$ permuting $p_1,...,p_r$ induces
a transitive action on the components of $\Hs$). Hence this image,
i.e., the locus of genus $g$ curves admitting a cover to $\P$ of type $\s$,
is irreducible.

\begin{Definition}(a) The moduli dimension of $\s$, denoted by $\md(\s)$,
is the dimension of the
image of $\Phi_\s$; i.e., the dimension of the locus of genus $g$
curves admitting a cover to $\P$ of type $\s$. We say $\s$ has
{\bf full moduli dimension} if\  $\md(\s)\ =\ \dim \M_g$.

(b) We say $\s$ has {\bf infinite moduli degree} if the following
holds: If $f:X\to\P$ is a cover of type $\s$ with general branch
points then $X$ has infinitely many covers to $\P$ of (the same)
type $\s$ such that the corresponding subfields of the function
field of $X$ are all different. (This terminology is further
discussed at the end of this section).
\end{Definition}

A curve is called a {\bf general curve of genus $g$} if it
corresponds to a point of $\M_g$ that does not lie in any proper
closed subvariety of $\M_g$ defined over $\bar\bQ$ (the algebraic
closure of the rationals). Clearly, an admissible tuple $\s$ has
full moduli dimension if and only if each  general curve of genus
$g_\s$ admits a cover to $\P$ of type $\s$.

Part (a) of the following Remark
is the necessary condition for full moduli dimension used by
Guralnick, Fried and Zariski. We indicate the proof at the end of
this section.

\begin{Remark} \label{nec} Let $\s$ be an admissible tuple of length
$r$ in $S_n$, and $g:=g_\s$.\newline
(a) Suppose $\s$ has full moduli
dimension. Then $r-3\ge\ \dim \M_g$, thus if $g\ge2$
then $r \ge 3g$.\newline
(b) 
 If\ $r-3>\ \dim \M_g$ then $\s$ has
infinite moduli degree.
\end{Remark}

Here is a  simple but crucial lemma that allows us to make use of the
hypothesis of infinite moduli degree.

\begin{Lemma} \label{Lemma2} Suppose $f_i:X\to\P$ is an infinite collection
of covers such that the  corresponding subfields of the
function field of $X$ are all different.
Let $S$ be the set of $(x,y)\in X\times X$ with $f_i(x)=f_i(y)$ for some $i$.
Then $S$ is Zariski-dense in $X\times X$.
\end{Lemma}

\Proof Let $S_i$ be the curve on $X\times X$ consisting of all
$(x,y)$ with $f_i(x)=f_i(y)$. The set $S$ is the union of all $S_i$.
If $S$ is not Zariski-dense in $X\times X$ then it must be the
union of finitely many $S_i$; then the curves $S_i$ cannot be all
distinct. But if $S_i=S_j$ then the subfields of $\bC(X)$ corresponding
to $f_i$ and $f_j$ coincide. This contradicts the hypothesis.
\qed

Here is our sufficient condition for full moduli dimension.

\begin{Lemma} \label{Lemma1} Let $n\ge3$.
Given an admissible tuple $\s=\ (\s_1,...,\s_r)$ in $S_n$
with $g_\s>0$, define
$\ss =\ (\s_1,...,\s_{r+2})$, where either $$\s_{r+1} \ \ = \ \
\s_{r+2} \ \ = \ \ (1,2)(n,n+1)$$ is a double transposition or
$$\s_{r+1} \ \ = \ \ \s_{r+2}^{-1} \ \ = \ \ (n-1,n,n+1)$$ is a
3-cycle. Then $\ss$ is an admissible tuple in $S_{n+1}$ with
$g_{\ss}=g_\s+1$. If $\s$ has infinite moduli degree then
 $$\md(\ss)\ \ \ \ge \ \ \ \md(\s)\ +\ \cases{3 & \ \ if\ \ $g_\s>1$\cr 2 & \ \ if\
\ $g_\s=1$}$$
\end{Lemma}

\Proof Let $g:=g_\s$.  Then $g_{\ss}=g+1$ by Riemann-Hurwitz.
 Let $\Phi:=\Phi_{\ss}$ and $\H:=\Hss$. The map $\Phi$
 extends to $\bar\Phi: \bar\H  \to  \bar\M_{g+1}$,
where $\bar\M_{g+1}$ is the stable compactification of $\M_{g+1}$,
and $\bar\H$ is $\H$ plus that piece $\partial\H$ of the boundary
where the last two branch points come together (see \cite{We});
thus $\bar\H$ covers the set of $(p_1,...,p_{r+2})$ in
$(\P)^{r+2}$ with $p_i\ne p_j$ for $i\ne j$ unless $\{i,j\}=
\{r+1,r+2\}$, and $\partial\H$ is the inverse image of the subset
defined by the condition $p_{r+1}=p_{r+2}$.

If we coalesce the last two entries of $\ss$ we obtain $\s$, which
has  orbits of length $n$ and 1 on $\{1,...,n+1\}$. For a cover
$X_{g+1}\to \P$ of type $\ss$, this means the following:  When
coalescing the last two branch points, $X_{g+1}$ degenerates into
a nodal curve $\bar X$ with two components linked at one point
$P$.
One component is a non-singular curve covering $\P$ of degree 1.
The other component $\bar X_{g}$ is a singular curve  whose only
singularity is a node $N$.
Its normalization $X_{g}$ covers $\P$ of type $\s$. If
$\s_{r+1}=(1,2)(n,n+1)$ then $N$ corresponds to the cycle $(1,2)$
and $P$ to the cycle $(n,n+1)$. If $\s_{r+1}=(n-1,n,n+1)$ then
$N=P$.

The nodal curve $\bar X$ is stably equivalent to the stable curve $\bar X_{g}$,
and the latter constitutes the image in $\bar \M_{g+1}$ of the element of $\partial\H$
corresponding to $\bar X\to\P$ (see \cite{HM}, Th. 3.160). Thus the image
of $\partial\H$ in $\bar \M_{g+1}$ lies in the boundary component consisting of
irreducible curves with one node whose normalization has genus $g$. We can
identify this boundary component with $\M_{g,2}$ (= moduli space of genus $g$
curves with two unordered marked points). The two marked points correspond to
the node.
 Thus we have the commutative diagram
$$\matrix{\H_\s & \mathop{\longrightarrow}\limits^{\Phi_\s} & \M_{g} \cr
\Big\uparrow & & \Big\uparrow\cr
\partial\H & \mathop{\longrightarrow} & \M_{g,2} \cr
\Big\downarrow & & \Big\downarrow\cr
\bar\H & \mathop{\longrightarrow}\limits^{\bar\Phi}& \bar\M_{g+1}\cr}$$
where the vertical arrows on the lower level are inclusion. The map $\M_{g,2}\to \M_{g}$
is the natural projection (forgetting the marked points), and the map
$\partial\H\to \H_\s$ sends the point corresponding to the cover $\bar X\to\P$
to that corresponding to the cover
$X_{g}\to \P$ of type $\s$ (see the previous paragraph).

The image of $\bar\H$ in $\bar\M_{g+1}$ is irreducible (as remarked above).
Its intersection with the boundary of $\bar\M_{g+1}$ is a closed
proper subvariety, hence has codimension at least 1. This subvariety contains
the image of $\partial\H$, which we denote by $\mbox{Im}(\partial\H)$.
Thus $\md(\ss)\ge \ 1+ \dim \mbox{Im}(\partial\H)$.

The fiber $F$ in $\M_{g,2}$ of the point of $\M_{g}$ corresponding to $X_{g}$
can be identified with the set of unordered pairs $(x,y)$ of distinct points
of $X_g$, modulo Aut$(X_{g})$. The intersection $F_\s$
of this fiber with $\mbox{Im}(\partial\H)$
consists of those  $(x,y)$ such that there is a cover
$f:X_{g}\to \P$ of type $\s$ with $f(x)=f(y)$ and $f(x)$ not a branch point
of $f$. If $X_{g}$ is a general curve with
the property that it admits a cover to $\P$ of type $\s$, then by Lemma 
\ref{Lemma2}
and the hypothesis of infinite moduli degree,
 $F_\s$ is Zariski-dense in $F$. Since $F_\s$ is the general
fiber of the surjective map $\mbox{Im}(\partial\H)\to \Phi_\s(\H_\s)$,
it follows that
$\dim \mbox{Im}(\partial\H)= \ \dim F + \dim \Phi_\s(\H_\s)=\ \dim F +\md(\s)$.
This completes the proof.\qed

\medskip
\def\PGL2{\mbox{PGL}_2(\bC)}
Consider the natural action of $\PGL2$  on $\P$ (by fractional
linear transformations). It induces an action on $\Hs$, with
$\lambda\in \PGL2$ mapping $([f], (p_1,...,p_r))$ to
$([\lambda\circ f], (\lambda(p_1),..., \lambda(p_r)))$. The closed
subspace of $\Hs$ defined by the conditions $p_1=0$, $p_2=1$,
$p_3=\infty$ maps bijectively to the quotient $\Hs/\PGL2$. Hence
this quotient carries a natural structure of quasi-projective
variety, and  the map $\Phi_\s: \Hs\to\M_g$ induces a morphism
$\Hs/\PGL2\to\M_g$. (Clearly $\Phi_\s$ is constant on
$\PGL2$-orbits).

The dimension of (each component of) $\Hs/\PGL2$ is $r-3$. Thus if
$\Phi_\s$ is dominant then \ $r-3\ge\ \dim \M_g$. This proves 
Remark \ref{nec}(a).
If $r-3>\ \dim \M_g$ then the general fiber of the map
$\Hs/\PGL2\to\M_g$ is infinite. This proves Remark \ref{nec}(b)
(since two covers $f_1,f_2:X\to\P$
correspond to the same subfield of the function field of $X$ if
and only if $f_1$ is the composition of $f_2$ with an element of
$\PGL2$). 

For clarification, we now briefly discuss the general concept of
moduli degree. This will not be needed elsewhere in the paper. The
map $\Hs/\PGL2\to\M_g$ factorizes further over the action of $S_r$
permuting the branch points (i.e., one can drop the ordering of
the branch points. Actually, the version of the Hurwitz space
without ordering of the branch points is more natural, see
\cite{Buch}, Ch. 10, but for the purpose of this paper we need the
ordering). Anyway, the natural definition of the moduli degree of
$\s$ is as follows: The degree of the induced map from the
(irreducible) variety $\Hs/(\PGL2\times S_r)$ to $\M_g$. Thus the
moduli degree of $\s$ is the number of covers $f:X\to\P$ of type
$\s$ modulo $\PGL2$, where $X$ corresponds to a (fixed) general
point in the image of $\Phi_\s$.

\section{Covers with monodromy group $A_n$}

\def\DT{\mbox{DT}} \def\DTA{\mbox{DTA}}
We consider admissible tuples $\s=(\s_1,...,\s_r)$ in $S_n$
such that each $\s_i$ is a double transposition (resp., 3-cycle). Then
$r=n+g-1\ge n-1$, where $g:=g_\s$ (by Riemann-Hurwitz). Let DT$(n,g)$ 
(resp., TC$(n,g)$ ) be the set
of these tuples $\s$; and let DTA$(n,g)$ (resp., TCA$(n,g)$ )
be the subset consisting of those $\s$
that generate $A_n$ (the alternating group). 

\begin{Lemma} \label{Lemma3}  (i) For each $n\ge4$ (resp., $n\ge6$) the set DT$(n,0)$
(resp., DTA$(n,0)$ ) is non-empty. \newline
(ii) The set TCA$(n,0)$ is non-empty for each $n\ge3$. 
\end{Lemma}

\Proof (i) For $n=4$ take $\s$ to consist of all double transpositions in $A_4$.
For $n=5$ take $\s=(\s_1,...,\s_4)$ such that $\s_1\s_2$ ($= (\s_3\s_4)^{-1}$)
is a 5-cycle. For $n=6$ use GAP (or check otherwise).

Assume now
$\s$ is in DTA$(n,0)$, and $n\ge6$. We may assume $\s_r=(1,2)(3,4)$. Replacing $\s_r$
by the two elements $(1,2)(n,n+1)$ and $(3,4)(n,n+1)$ yields a
tuple in DTA$(n+1,0)$. This proves (i).

(ii) 
For $n=3$ take $\s=((1,2,3),(1,2,3)^{-1})$.
Assume now $\s$ is in TCA$(n,0)$, $n\ge3$. We may assume $\s_1=(1,2,3)$.
Replacing $\s_1$ by the two elements $(n+1,3,1)$ and $(3,n+1,2)$ yields a
tuple in TCA$(n+1,0)$. 

\begin{Lemma} \label{Lemma4} Both of DTA$(n,g)$ and TCA$(n,g)$ 
contain a tuple of full moduli dimension if one of the following
holds:\newline
(i) \ $g=1$ and $n\ge5$.\newline
(ii)\ $g=2$ and $n\ge6$.\newline
(iii)\ $g>2$ and $n\ge 2g+1$.
\end{Lemma}

\Proof (i) 
See \cite{FKK} for a proof of the TCA$(n,1)$ case that does not use
the stable compactification. For the DTA$(n,1)$ case, we use induction
on $n$. 

We anchor our induction at $n=5$. We choose
$\s = (\s_1,\s_2,\s_3,\s_4,\s_5)$, where $\s_1 = \s_2= (1,2)(3,4)$,
$\s_3 = (1,2)(4,5)$, $\s_4 = (1,4)(2,5)$, and $\s_5 = (1,5)(2,4)$. 
If we coalesce the last two entries of $\s$ we obtain $(\s_1,\s_1,\s_3,\s_3)$,
which
has orbits $\{1,2\}$ and $\{3,4,5\}$. For a cover
$X_{1}\to \P$ of type $\s$, this means the following:  When
coalescing the last two branch points, $X_{1}$ degenerates into
a nodal curve $\bar X$ with two components linked at one point
$P$.
Both components are non-singular curves of genus 1 (resp., 0). They 
both cover $\P$ with four branch points and of degree 2 (resp. 3).
The point $P$ ramifies in both covers. 
The nodal curve $\bar X$ is stably equivalent to its genus 1 component,
and the latter constitutes the image in $\M_{1}$ of the cover 
$\bar X\to\P$ (as in the proof of lemma \ref{Lemma1}). Clearly, 
every element of $\M_{1}$ can be obtained in this fashion. Thus 
the map $\H_{\s} \to \M_1$ is dominant because the boundary 
of $\H_{\s}$ already maps surjectively to $\M_{1}$.   


Now assume $\s=(\s_1,...,\s_{n})$ is a tuple in DTA$(n,1)$, $n\ge5$,
of full moduli dimension. Write $\s_{n}=st$ where $s,t$ are double
transpositions in $S_{n+1}\setminus S_n$. Let 
$\s' =\ (\s_1,...,\s_{n-1},s,t)$, a tuple in DTA$(n+1,1)$.
Moreover, $\s'$ has full moduli dimension because ${\Phi}_{\s'}$
restricted to the boundary component of ${\cal H}_{\s'}$ isomorphic
to ${\cal H}_\s$ already  maps dominantly to $\M_1$.

\smallskip\noindent
(ii) Same for both cases. So we only do the DT case.
By (i), there is a  tuple in  DTA$(n-1,1)$ of full moduli dimension.
Its length equals $n-1$, and $(n-1)-3>1=\dim \M_1$; thus the
tuple has infinite moduli degree by Remark \ref{nec}(b). Then
Lemma \ref{Lemma1} produces a tuple in DTA$(n,2)$ of 
full moduli dimension.

\smallskip\noindent
(iii)  Same for both cases. So we only do the DT case.
First we settle the case $g=3$, $n\ge 7$. By (ii),
there is a tuple in DTA$(n-1,2)$ and of moduli dimension 3.
Its length is $n$, and $n-3\ >3=\ \dim \M_2$; the claim follows
from Remark  \ref{nec}(b) and Lemma \ref{Lemma1}.

Now suppose $g>3$, $n\ge 2g+1$. Then $n-1\ge 2(g-1)+2$. By
induction we may assume there is a tuple in DTA$(n-1,g-1)$ and of
full moduli dimension. Its length is $r:=n+g-3$, and $r-3>\
3(g-1)-3= \ \dim \M_{g-1}$; the claim follows again from Remark
\ref{nec}(b) and Lemma \ref{Lemma1}.

\begin{Theorem} \label{thm}
 (i) Let $g\ge3$. Then each general curve of genus $g$ admits a cover to $\P$
of degree $n$ with monodromy group $A_n$ such that all inertia
groups are generated by double transpositions 
 if and only if $n\ge
2g+1$.
\newline (ii) For $n\ge 6$ (resp., $n\ge 5$), each general
curve of genus 2 (resp., 1) admits a cover to $\P$ of degree $n$
with monodromy group $A_n$ such that all inertia groups are
generated by double transpositions.
\newline (iii) Assertions (i) and (ii) also hold for 3-cycles
instead of double transpositions.
\end{Theorem}

\Proof In view of Lemma \ref{Lemma4}, it
only remains to show that the condition $n\ge 2g+1$ in (i) is
necessary. Indeed, if the general curve of genus $g$ admits such a
cover then an associated tuple of branch cycles is in DTA$(n,g)$
and of full moduli dimension. Thus the claim follows from the
necessary condition $r\ge3g$ (Remark \ref{nec}) since $r=n+g-1$.
The proof of (iii) is the same.

\begin{Corollary} Let $C$ be a general curve of genus $g\ge4$.
Then the monodromy groups of primitive covers $C\to\P$ are among the symmetric and
alternating groups, and up to finitely many, all of these groups occur.
\end{Corollary}

Here a cover is called primitive if it does not factor non-trivially.
The first assertion in the Corollary follows from \cite{GM},
and the second from the Theorem plus Brill-Noether theory.

\section{Braid orbits of  admissible tuples}

The  {\bf braid orbit} of a tuple $\s$ in $S_n$
is the smallest set of
tuples in $S_n$ that contains $\s$ and is closed under
(component-wise) conjugation and  under
the braid operations
$$ (g_1,...,g_r)^{Q_i}\ \ = \ \
(g_1,\ldots,\ g_{i+1},\ g_{i+1}^{-1}
 g_i g_{i+1} \ ,\ldots, g_r) $$
for $i=1,...,r-1$.

Let $\s$, $\s'$ be admissible tuples in $S_n$ of length $r$.
Let $f:X\to\P$ be a cover of type $\s$. Then  $f$ is of type $\s'$
if and only if $\s'$ lies in the braid orbit of $\s$.
In other words, for the associated Hurwitz spaces we have
${\cal H}_\s={\cal H}_{\s'}$ if and only if $\s'$ lies in the 
braid orbit of $\s$ (see \cite{FrV},  \cite{Buch}, Ch. 10).
Thus the above notions of  moduli dimension,  moduli degree etc.
depend only on the  braid orbit of $\s$. So from now on we will
speak of the  moduli dimension of a  braid orbit, etc.

\subsection{Braid orbits of 2-cycle tuples}

Admissible tuples in $S_n$ of fixed length  that consist only of
transpositions form a single braid orbit (by Clebsch 1872, see \cite{Buch},
Lemma 10.15). They correspond to the so-called {\bf simple covers}.
Their  braid orbit has full  moduli dimension if and only if $2(n-1)\ge g$,
where $g=g_\s$ (see the remarks in the Introduction).

\subsection{Braid orbits  of  3-cycle tuples}

Now consider tuples that consist only of 3-cycles.
Recall our notation TC$(n,g)$ for the set of those 
(admissible) tuples
with fixed parameters $n$, $g$. Assume $n\ge 5$.
Note that  TC$(n,g)=$ TCA$(n,g)$ (i.e., each such tuple generates
$A_n$) by \cite{Hup}, Satz 4.5.c and the fact that a transitive group 
generated by $3$-cycles must be primitive.
 The corresponding covers have
been studied by Fried \cite{Fr1}. Serre  \cite{Se1}, \cite{Se2}
considered certain generalizations. Fried proved that  TC$(n,g)$
(is non-empty and)
consists of exactly two braid orbits (resp., one  braid orbit)
if $g>0$ (resp., $g=0$). 
Let 
$$\{\pm1\}\ \ \to\ \ \hat A_n\ \ \to A_n\ \ $$ be the unique
non-split degree 2 extension of $A_n$. Each 3-cycle $t\in A_n$ 
 has a  unique lift $\hat t\in \hat A_n$ of order 3.
For $\s=(\s_1,...,\s_r)\in $ TC$(n,g)$ we have
$\hat \s_1\cdots \hat \s_r=\pm1$. The value of this product
is called the  {\bf lifting invariant} of $\s$. It depends only
on the braid orbit of $\s$. For $g=0$  the  lifting invariant
is $+1$ if and only if $n$ is odd (by \cite{Fr1} and \cite{Se1}).
For $g>0$ the two braid orbits on TC$(n,g)$ 
have distinct lifting invariant.

Now we can refine Theorem \ref{thm} as follows.

\begin{Theorem} \label{thm3}
Assume $n\ge 6$, $g>0$ and $n\ge2g+1$.
Then both  braid orbits on TC$(n,g)$
have full moduli dimension.
\end{Theorem}

\Proof The claim holds for $g=1$ by \cite{FKK},
Comment 0. Now
 suppose in the situation of Lemma \ref{Lemma1},
$\ss$ is a tuple in $A_n$ with
$\s_{r+1}  = \s_{r+2}^{-1} =  (n-1,n,n+1)$. Then clearly
$\s$ and $\tilde\s$ have  the same lifting invariant.
Thus the proof of Lemma \ref{Lemma4} also shows the
present refinement, since it iterates the construction
of Lemma \ref{Lemma1}.


\section{The exceptional cases in genus 3}
\def\t{\tau}

Let $\s=(\s_1,...,\s_r)$ be an admissible tuple in $S_n$, and
$g:=g_\s\ge3$. Assume $\s$ satisfies the necessary condition $r\ge
3g$ for full moduli dimension. Assume further $\s$ generates a
primitive subgroup $G$ of $S_n$. If $g\ge4$ then $G=S_n$ or
$G=A_n$ by \cite{GM} and \cite{GS}. If $g=3$ and $G$ is not 
$S_n$ or $A_n$ then  one of the following holds  (see
\cite{GM}, Theorem 2):

\begin{enumerate}
\item[(1)]  \ \ $n=7$, \ $G\cong GL_3(2)$
\item[(2)]  \ \  $n=8$, \ $G\cong AGL_3(2)$ (the affine group)
\item[(3)] \ \   $n=16$, \ $G\cong AGL_4(2)$
\end{enumerate}

\noindent  Recall that $GL_3(2)$ is a simple group of order 168.
It acts doubly transitively on the 7 non-zero elements of
$(\bF_2)^3$. The affine group $ AGL_m(2)$ is the semi-direct product
of $ GL_m(2)$ with the group of translations; it acts triply
transitively on the affine space $(\bF_2)^m$.

In cases (1) and (3), the tuple $\s$ consists of 9
 transvections of the respective linear or affine
group. (A transvection fixes a hyperplane of the underlying linear
or affine space point-wise). In case (2), either $\s$ consists of
10 transvections or it consists of 8 transvections plus an element
of order 2,3 or 4 (where the element of order 2 is a translation).

\begin{Remark}  The tuples in case (1) form a single braid orbit
on DT$(7,3)$. This  braid orbit has full moduli dimension by the
Theorem below.
\end{Remark}

\Proof We show that tuples of 9 involutions generating $G=GL_3(2)$
(with product 1) form a single braid orbit. This uses the BRAID program
\cite{MS}. 
Direct application of the program is not
possible because the number of  tuples is too large.

We first note that if  9 involutions generate $G$, then there are
6 among them that  generate already (since the lattice of subgroups
of $G$ has length 6). We can move these 6 into the first 6 positions
of the tuple by a sequence of braids. Now we apply the BRAID program
to 6-tuples of involutions generating $G$ (but not necessarily with
product 1). We find that such tuples with any prescribed value of their
product form  a single braid orbit. By inspection of these braid orbits,
we find that each contains a tuple whose first two involutions are
equal, and the remaining still generate $G$. This reduces the original
problem to showing that  tuples of 7 involutions with product 1,
generating $G$, form a single braid orbit. The BRAID program did that.
\qed

\medskip
In cases (1) and (2), the transvections yield double
transpositions in $S_n$. Thus again Lemma \ref{Lemma1} can be used
 to show there actually exist such tuples that have full moduli
dimension. Case (3) requires a more complicated argument which
will be worked out later.

\begin{Theorem} \label{thm2} Each general curve of
genus 3 admits a cover to $\P$ of degree 7 (resp., 8) and
monodromy group $GL_3(2)$ (resp., $ AGL_3(2)$ ), branched at 9
(resp., 10) points of $\P$, such that all inertia groups are generated by
double transpositions.
\end{Theorem}

\Proof  Let $G$ be a (doubly) transitive subgroup of $S_7$
isomorphic to $GL_3(2)$. Let $H$ ($\cong S_4$) be a point
stabilizer in $G$. View $H$ as a subgroup of $S_6$ via its
(transitive) action on the other 6 points. In \ref{genus2} below,
we show there is a tuple $\t$ in $\DT(6,2)$ of full moduli
dimension that generates this subgroup $H$ of $S_6$. This tuple
has length 7, hence has infinite moduli degree by Remark
\ref{nec}(b). Choose a double transposition in $G$ that is not in
$H$, and append two copies of it to the tuple $\t$. By Lemma
\ref{Lemma1}, this yields a tuple $\s\in\DT(7,3)$ of full moduli
dimension, satisfying (1).

\smallskip\noindent
The group $GL_3(2)$ is the stabilizer of 0 in the transitive
action of $ AGL_3(2)$ on the 8 points of $(\bF_2)^3$. Replacing the
last entry $\s_9$ of the above tuple $\s$ 
by two double
transpositions from $ AGL_3(2)$ that are not in $ GL_3(2)$ and
have product $\s_9$, yields a tuple in $\DT(8,3)$ satisfying (2).
This tuple has full moduli dimension because already the boundary
of the corresponding Hurwitz space maps dominantly to $\M_3$.




\subsection{Certain covers of degree 6 from the general curve of genus 2 to $\P$}
\label{genus2}

Let $\tau_1,\tau_2,\tau_3$ be the three double transpositions in
$H:=S_4$. Let $\rho_1$ and $\rho_2$ be transpositions in $H$
generating an $S_3$-subgroup. Then the tuple $$\tau \ \ =\ \
(\tau_1,\tau_2,\tau_3,\rho_1,\rho_1,\rho_2,\rho_2)$$  generates
$H$. View $H$ as a subgroup of $S_6$ as in the proof of Theorem
\ref{thm2}. Then $\tau$ becomes an element of $\DT(6,2)$ (since
all involutions of $GL_3(2)$ act as double transpositions on the 7
points).

Now consider a cover $f:X\to\P$ of type $\tau$. Note that $H$ is
an imprimitive subgroup of $S_6$, permuting 3 blocks of size 2.
The kernel of the action of $H$ on these 3 blocks equals
$\{1,\tau_1,\tau_2,\tau_3\}$. Thus $f$ factors as $f=hg$ where
$g:X\to\P$ is of degree 2 (the hyperelliptic map on the genus 2
curve $X$) and $h:\P\to\P$ is a simple cover of degree 3 (i.e.,
its tuple of branch cycles consists of 4 involutions in $S_3$).
Let $p_i\in\P$, $i=1,2,3$ be the branch point of $f$ corresponding
to $\tau_i$. Then $p_i$ has 3 distinct pre-images $x_i,y_i,z_i$
under $h$. We may assume $p_1=0=x_3$, $p_2=\infty=y_3$,
$p_3=1=x_1$.  Then $h$ is of the form $$h(x)\ \ = \ \
\frac{(x-1)(x-y_1)(x-z_1)}{(x-x_2)(x-y_2)(x-z_2)}$$ Exactly one
of $x_i,y_i,z_i$, say $z_i$, is unramified under $g$. Thus
$x_1=1,y_1,x_2,y_2,x_3=0,y_3=\infty$ are the 6 branch points of the
hyperelliptic map $g$. It is well-known that (the $\PGL2$-orbit
of) this 6-set determines the isomorphism class of the genus 2
curve $X$. Now we are ready to prove:

\begin{Lemma} \label{g2}
The tuple $\tau$  has full moduli dimension.
\end{Lemma}

\Proof It
suffices to show that for each choice of $y_1',x_2',y_2'$
sufficiently close to $y_1,x_2,y_2$, respectively (in the
complex topology), the following holds: There are $z_1',z_2'$
close to $z_1,z_2$, respectively, such that the map $$h'(x)\ \ = \
\ \frac{(x-1)(x-y'_1)(x-z'_1)}{(x-x'_2)(x-y'_2)(x-z'_2)}$$
composed with the double cover $g':X'\to\P$ branched at
$y_1',x_2',y_2',0,\infty,1$ is a cover of type $\t$. This
follows by continuity once we know that the condition $h'(0)=1 \
(= h'(\infty) )$ is preserved. But this condition $h'(0)=1$ is
easy to achieve: We can view it as defining $z_2'$ (after free
choice of $z_1'$).


\end{document}